\documentclass{article}
\usepackage{amssymb,amsmath,amsthm,graphicx,xcolor,tikz}
\usetikzlibrary{knots}
\usetikzlibrary{hobby}
\usetikzlibrary{arrows.meta}
\usetikzlibrary{arrows,decorations.markings}
\def\utr{\, \underline{\triangleright}\, }
\def\otr{\, \overline{\triangleright}\, }
\def\ud{\, \underline{\bullet}\, }
\def\od{\, \overline{\bullet}\, }

\newtheorem{theorem}{Theorem}

\newtheorem{proposition}[theorem]{Proposition}

\theoremstyle{definition}
\newtheorem{example}{Example}
\newtheorem{definition}{Definition}
\newtheorem{remark}{Remark}

\date{}

\title{\Large \textbf{Cocycle Enhancements of Psyquandle Counting Invariants}}

\author{Jose Ceniceros\footnote{Email: jcenicer@hamilton.edu}\and
Sam Nelson\footnote{Email: Sam.Nelson@cmc.edu. 
Partially supported by Simons Foundation collaboration grant 702597.}}

\begin{document}
\maketitle

\begin{abstract}
We bring cocycle enhancement theory to the case of psyquandles. Analogously
to our previous work on virtual biquandle cocycle enhancements, we define
enhancements of the psyquandle counting invariant via pairs of a biquandle
2-cocycle and a new function satisfying some conditions. As an application we 
define new single-variable and two-variable polynomial invariants of
oriented pseudoknots and singular knots and links. We provide examples
to show that the new invariants are proper enhancements of the counting 
invariant are are not determined by the Jablan polynomial.
\end{abstract}

\parbox{5.5in} {\textsc{Keywords:} Psyquandles, cocycle invariants, pseduoknots,
singular knots

                \smallskip
                
                \textsc{2020 MSC:} 57K12}

\section{Introduction}

In \cite{CJKLS}, a (co)-homology theory for quandles was introduced and
used to enhance the quandle counting invariant of classical (and later 
virtual) knots and links with $2$-cocycles. The idea was extended to 
biquandles in \cite{CES} and to virtual biquandles by the present authors
in \cite{CN1}, obtaining a two-variable enhancement in the case of
\textit{strongly compatible} cocycles and a single-variable polynomial
enhancement in the case of \textit{weakly compatible} cocycles.

The term \textit{pseudoknot} was first used in the natural sciences to describe
knotted structures with only partial crossing information, such as DNA 
strands in images with insufficient resolution to determine the crossing
information at some crossings, for example \cite{CDDK}. A rigorous 
mathematical definition
for pseudoknots was established in \cite{H} and has been studied in further 
work such as \cite{HHJJMR, HJ, HJJ}. In addition to classical crossings, 
pseudoknots 
have \textit{precrossings} which are understood to be classical crossings
about which we lack crossing information. In particular, a pseduoknot can 
be understood as a kind of probability distribution with the classical knot
resolutions obtained by assigning crossing information to precrossings
as outcomes.

\textit{Singular knots} are rigid vertex isotopy classes of 4-valent spatial
graphs. They have been studied particularly for their connection with
Vassiliev invariants; see \cite{V2} for example. In particular, the singular
Reidemeister moves are the same as the pseduoknot Reidemeister moves, aside
from one move, if we replace precrossings with singular crossings with
precrossings.

In \cite{NOS}, an algebraic structure called \textit{psyquandle} was introduced,
algebraically encoding the Reidemeister moves for pseudoknots and singular
knots into one unified structure. Psyquandles were used to define a new 
polynomial invariant of pseudoknots and singular, the Jablan polynomial, which
can be understood as a weighted sum of Alexander polynomials of the 
classical resolutions of the pseudoknot. An integer-valued psyquandle 
counting invariant was also introduced.

In this paper we enhance the psyquandle counting invariant with cocycles 
analogously to our work in \cite{CN1}, obtaining a new infinite family of 
single-variable and two-variable polynomial invariants of pseudoknots
and singular knots and links. The paper is organized as follows. In Section 
\ref{P} we revisit psyquandles and recall their basics. In Section \ref{BC}
we recall the basics of biquandle (co)homology. In Section \ref{CE} we
introduce the new invariants and provide examples and computations, including
and example to show that the new invariants are not determined by the Jablan
polynomial. We conclude in Section \ref{Q} with some questions for future 
research.

\section{Psyquandles}\label{P}

In this section we recall the basics of psyquandles. See \cite{NOS} for more.

\begin{definition}
Let $X$ be a set. A \textit{psyquandle structure} on $X$ is set of four
binary operations $\utr,\otr,\ud,\od:X\times X\to X$ satisfying the conditions
\begin{itemize}
\item[(0)] All four operations are right-invertible,i.e. there exist
binary operations $\utr^{-1},\otr^{-1},\ud^{-1},\od^{-1}:X\times X\to X$ such that
\[\begin{array}{rcccl}
(x\utr y)\utr^{-1}y & = & (x\utr^{-1}y)\utr y & = & x\\ 
(x\otr y)\otr^{-1}y & = & (x\otr^{-1}y)\otr y & = & x\\ 
(x\ud y)\ud^{-1}y & = & (x\ud^{-1}y)\ud y & = & x\\ 
(x\od y)\od^{-1}y & = & (x\od^{-1}y)\od y & = & x,
\end{array}\]
\item[(i)] For all $x\in X$, $x\utr y=x\otr y$, 
\item[(ii)] For all $x,y\in X$, the maps $S,S':X\times X\to X\times X$ 
defined by
\[S(x,y)=(y\otr x,x\utr y) \quad\mathrm{and}\quad S'(x,y)=(y\od x,x\ud y)\] 
are invertible,
\item[(iii)] For all $x,y,z\in X$,
\[\begin{array}{rcl}
(x\utr y)\utr (z\utr y) & = & (x\utr z)\utr (y\otr z) \\
(x\utr y)\otr (z\utr y) & = & (x\otr z)\utr (y\otr z) \\
(x\otr y)\otr (z\otr y) & = & (x\otr z)\otr (y\utr z) \\
\end{array}\]
\item[(iv)] For all $x,y\in X$ we have
\[\begin{array}{rcl}
x\ud((y\otr x)\od^{-1} x) & = & [(x\utr y)\od^{-1} y]\otr[(y\otr x)\ud^{-1} x]\\
y\ud((x\utr y)\od^{-1} y) & = & [(y\otr x)\od^{-1} x]\utr[(x\utr y)\od^{-1} y],
\end{array}\]
and
\item[(v)] For all $x,y,z\in X$ we have
\[\begin{array}{rcl}
(x\otr y)\otr (z\od y) & = & (x\otr z)\otr (y\ud z) \\
(x\utr y)\utr (z\od y) & = & (x\utr z)\utr (y\ud z) \\
(x\otr y)\od (z\otr y) & = & (x\od z)\otr (y\utr z) \\
(x\utr y)\ud (z\utr y) & = & (x\ud z)\utr (y\otr z) \\
(x\otr y)\ud (z\otr y) & = & (x\ud z)\otr (y\utr z) \\
(x\utr y)\od (z\utr y) & = & (x\od z)\utr (y\otr z).
\end{array}\]
\end{itemize}
A psyquandle which also satisfies $x\ud x=x\od x$ for all
$x\in X$ is said to be \textit{pI-adequate}.
\end{definition}

\begin{example}
Every biquandle is a psyquandle by setting $x\od y=x\otr y$ and
$x\ud y=x\utr y$. 
\end{example}

\begin{example}
A module over $\mathbb{Z}[t^{\pm 1},s^{\pm 1},a^{\pm 1},b^{\pm 1}]/(t+s-a-b)$ 
is a psyquandle with operations
\[\begin{array}{rcl}
x\utr y & = & tx+(s-t)y \\
x\otr y & = & sx \\
x\ud y & = & ax+(s-a) y \\
x\od y & = & bx+(s-b) y
\end{array}\]
known as an \textit{Alexander psyquandle}.
\end{example}

\begin{example}
Given a finite set $X=\{1,2,\dots, n\}$, we can specify a psyquandle
structure on $X$ by explicitly listing the operation tables of the four
psyquandle operations. In practice it is convenient to put these together 
into an $n\times 4n$ bock matrix, so the psyquandle structure on $X=\{1,2,3\}$
specified by
\[
\begin{array}{r|rrr} \utr & 1 & 2 & 3 \\ \hline 1 & 2 & 2 & 2 \\ 2 & 3 & 3 & 3 \\ 3 & 1 & 1 & 1\end{array}
\begin{array}{r|rrr} \otr & 1 & 2 & 3 \\ \hline 1 & 2 & 2 & 2 \\ 2 & 3 & 3 & 3 \\ 3 & 1 & 1 & 1\end{array}
\begin{array}{r|rrr} \ud & 1 & 2 & 3 \\ \hline 1 & 3 & 3 & 3 \\ 2 & 1 & 1 & 1 \\ 3 & 2 & 2 & 2\end{array}
\begin{array}{r|rrr} \od & 1 & 2 & 3 \\ \hline 1 & 3 & 3 & 3 \\ 2 & 1 & 1 & 1 \\ 3 & 2 & 2 & 2\end{array}
\]
is encoded as the block matrix
\[
\left[\begin{array}{rrr|rrr|rrr|rrr}
2 & 2 & 2 & 2 & 2 & 2 & 3 & 3 & 3 & 3 & 3 & 3 \\
3 & 3 & 3 & 3 & 3 & 3 & 1 & 1 & 1 & 1 & 1 & 1 \\
1 & 1 & 1 & 1 & 1 & 1 & 2 & 2 & 2 & 2 & 2 & 2
\end{array}\right].
\]
\end{example}

See \cite{NOS} for more examples of psyquandles.

The psyquandle axioms are motivated by the semiarc coloring rules for
singular knots and pseudoknots below:
\begin{equation}\includegraphics{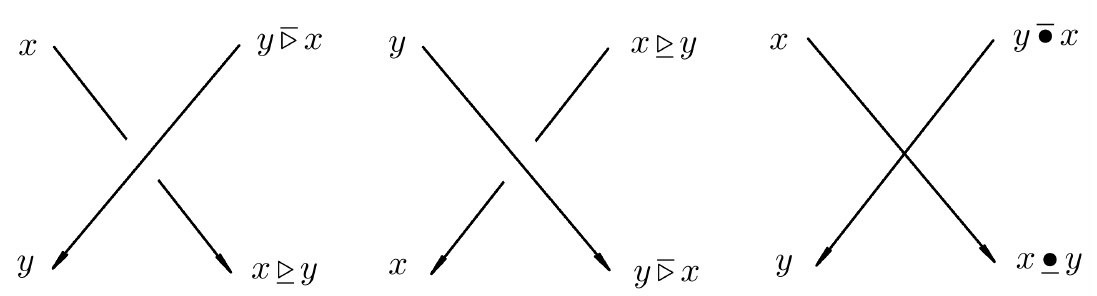}\label{eq1}\end{equation}
\begin{equation}\includegraphics{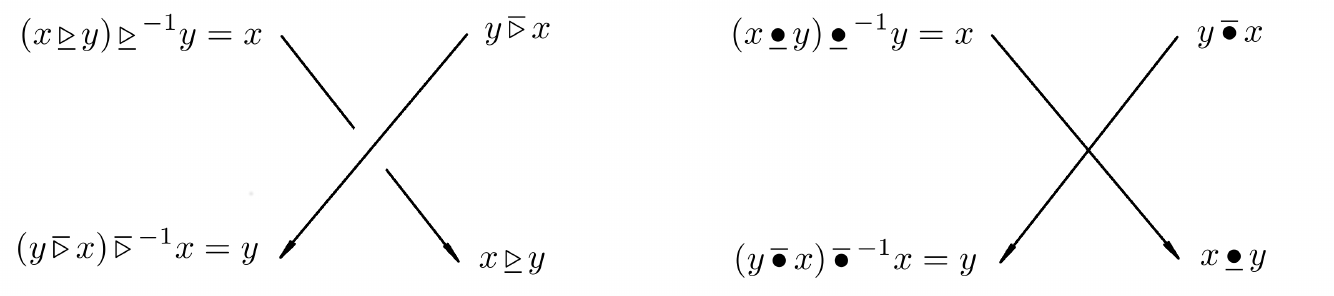}\label{eq2}\end{equation}

\begin{remark}
We have reformulated axiom (iv) from its description in \cite{NOR} in order 
to simplify the Reidemeister IV conditions for the Boltzmann weights we 
will define later. To see how these arise, consider the labeled move below:
\[\scalebox{0.7}{\includegraphics{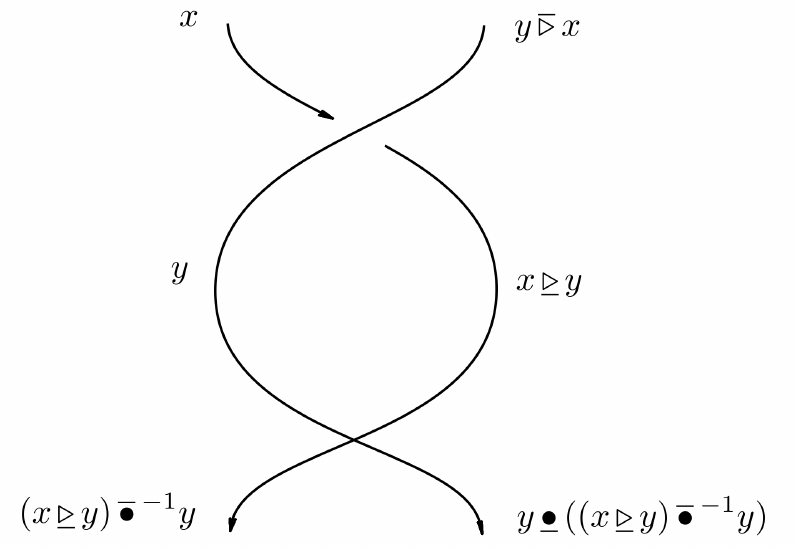}}
\quad \raisebox{0.7in}{$\leftrightarrow$}\quad
\raisebox{-0.4in}{\scalebox{0.7}{\includegraphics{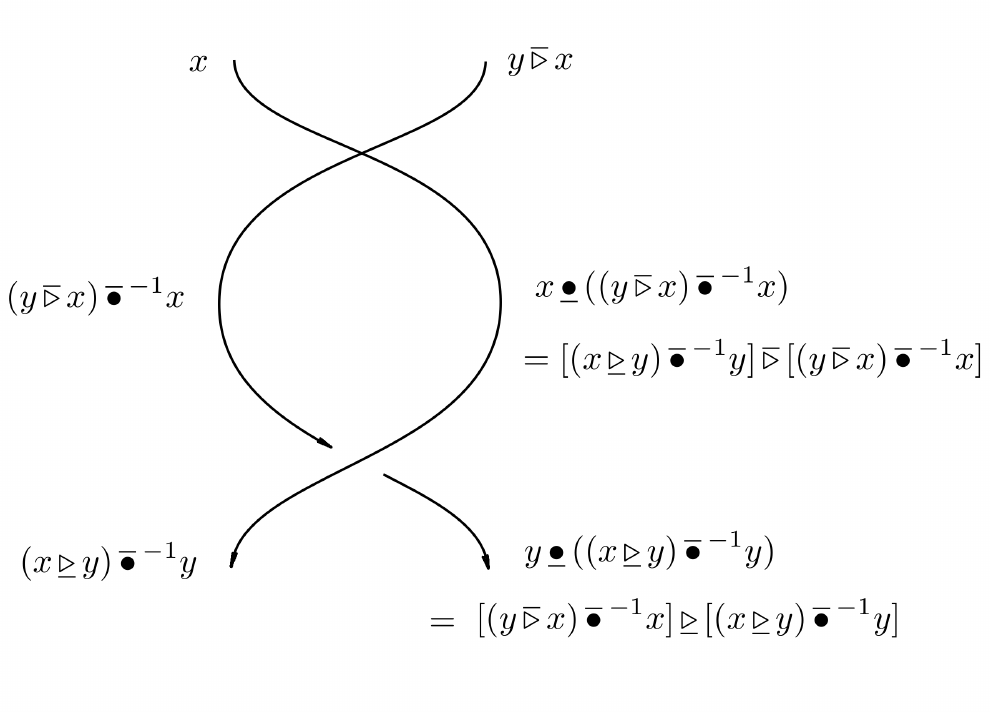}}}.
\]
\end{remark}

\begin{definition}
Let $X$ be a finite psyquandle (respectively, a $pI$-adequate finite 
psyquandle). An \textit{$X$-coloring} of an oriented singular 
knot or link diagram $L$ (respectively, an oriented pseudoknot or 
pseudolink diagram $L$) is an assignment of an element of $X$ to each semiarc
in $L$ such that the coloring rules (\ref{eq1}) and (\ref{eq2}) are
satisfied at every crossing.
\end{definition}

In \cite{NOS} we find the following result:
\begin{theorem} 
Let $X$ be a finite psyquandle (respectively, a finite $pI$-adequate psyquandle)
and let $L$ be an oriented singular knot or link diagram (respectively, an 
oriented pseudoknot or pseudolink diagram). Then the cardinality of the set of
$X$-colorings of $L$, 
\[\Phi_X^{\mathbb{Z}}(L)=|\mathrm{Hom}(P(L),X)|\]
is an integer-valued invariant of singular knots and links (respectively, 
pseudoknots and pseudolinks).
\end{theorem}

In the remainder of this paper we will define Boltzmann weight enhancements
of the psyquandle counting invariant.

\section{Biquandle Cohomology}\label{BC}

In this section we recall the basics of biquandle homology and cohomology.
See \cite{CES, EN} etc. for more on this topic.

Let $X$ be a finite biquandle and $R$ a commutative ring with identity. 
The set $C_n(X;R)=R[X^n]$ is the free $R$-module on the set ordered $n$-tuples
of element of $X$; its dual, $C_n(X;R)=\mathrm{Hom}(C_n(X), R)$ is the
set of $R$-module homomorphisms from $C_n(X;R)$ to $R$.

The map $\partial_n:C_n(X;R)\to C_{n-1}(X;R)$ defined on generators by 
\[\partial(x_1,\dots, x_n)=\sum_{k=1}^n (-1)^k(
(x_1,\dots, x_{k-1},x_{k+1,\dots, x_n})
-(x_1\utr x_k,\dots, x_{k-1}\utr x_k, x_{k+1}\otr x_k, \dots, x_n\otr x_k))\]
and extended linearly is a boundary map, with coboundary map
$\delta^n:C^n(X;R)\to C^{n-1}(X;R)$ given by $\delta^n(f)=f\partial_n$.
The resulting homology groups
$H_n=\mathrm{Ker}\partial_3/\mathrm{Im} \partial_2$ and
cohomology groups $H^n=\mathrm{Ker}\delta^3/\mathrm{Im} \delta^2$ are known as 
the \textit{birack homology} and \textit{birack cohomology} groups of $X$
with coefficients in $\mathbb{R}$. 

The subgroups $D_n(X;R)$ and $D^n(X;R)$ generated by elements 
$(x_1,\dots, x_n)$ with $x_k=x_{k+1}$ for some $k\in{1,\dots, n-1}$ form
\textit{degenerate} subcomplexes; modding out by these subcomplex yields 
\textit{biquandle homology} and \textit{biquandle cohomology}.

\begin{example}
A function $\phi:X\times X\to R$ represents a biquandle cohomology class 
in the biquandle cohomology of $X$ with $R$ coefficients if it satisfies
\begin{itemize}
\item[(i)] For all $x\in X$,
\[\phi(x,x)=0\]and
\item[(ii)] For all $x,y,z\in X$,
\[
\phi(x,y)
-\phi(x\utr z,y\utr z)
-\phi(x,z) 
+\phi(x\utr y,z\otr y)
+\phi(y,z)
-\phi(y\otr x,z\otr x)=0.\]
\end{itemize}
\end{example}

Biquandle 2-cocycles are of interest since they can be used to \textit{enhance}
the biquandle counting invariant, resulting in an invariant of oriented knots 
and links known as a \textit{cocycle enhancement}. More precisely, consider the
set of biquandle colorings of an oriented link diagram (i.e., satisfying 
equation (\ref{eq1}) at every crossing). At each crossing, we collect a 
contribution of $\pm \phi(x,y)$, with the resulting sum known as the 
\textit{Boltzmann weight} of the coloring. It is straightforward to check that
the biquandle cocycle conditions imply that such a Boltzmann weight is not
changed by Reidemeister moves; hence, the multiset of Boltzmann weights
forms an enhanced invariant of oriented links. It is common to encode these
multisets as ``polynomials'' by making multiset elements powers of a formal
variable with multiplicities as coefficients, for ease of comparison.

A biquandle 2-cocycle for a finite biquandle structure on $X=\{1,\dots, n\}$ 
can be written as a linear combination of characteristic functions 
\[\phi=\sum_{(j,k)\in X\times X}\phi_{jk}\chi_{(j,k)}\]
which we can conveniently encode as an $n\times n$ matrix whose
$(j,k)$ entry is $\phi_{jk}$.

\section{Cocycle Enhancements}\label{CE}

We will now generalize biquandle cocycle invariants to the case of psyquandles.

\begin{definition}\label{CocyclesDef}
Let $X$ be a psyquandle and $R$ a commutative ring with identity. A 
\textit{Boltzmann weight} for $X$ is a pair of maps $\phi,\psi:X\times X\to R$
satisfying
\begin{itemize}
\item[(i)] For all $x\in X$, $\phi(x,x)=0$
\item[(ii)] For all $x,y\in X$,
\[\phi(x,y)+\psi(y,(x\utr y)\od^{-1}y)=\phi((y\otr x)\od^{-1} x, (x\utr y)\od^{-1}y)+\psi(x,(y\otr x)\od^{-1} x).\]
\item[(iii)] For all $x,y, z \in X$,
\begin{eqnarray*}
\phi(x,y)+\phi(y,z)+\phi(x\utr y,z\otr y) & = & \phi(x\utr z,y\utr z)+\phi(x,z) + \phi(y\otr x,z\otr x)\\
\psi(x,y)+\phi(y,z)+\phi(x \ud y, z\otr y) &=& \psi(x \utr z, y \utr z)+\phi(x,z) + \phi(y \od x, z \otr x) \\
\psi (z,y) - \phi(x,y) - \phi(x \utr y, z \ud y) &=& \psi(z \otr x, y \otr x) -\phi(x,z) - \phi(x \utr z, y \od z).
\end{eqnarray*}
\end{itemize}
If $X$ is pI-adequate, we say that $(\phi,\psi)$ is \textit{$pI$-adequate} if 
\begin{itemize}
\item[(v)] For all $x \in X$,
\[\psi(x,x) = 0\]
\end{itemize}
and we say that $\phi$ and $\psi$ are \textit{strongly compatible} if
we also have
\begin{itemize}
\item[(vi)] For all $x,y,z\in X$,
\[\psi(x,y)=\psi(x\utr z, y\utr z)\quad\mathrm{and}\quad 
\psi(z,y)=\psi(z\otr x, y\otr x).\]
\end{itemize}
\end{definition}

The Boltzmann weight axioms are motivated by the Reidemeister moves for 
singular knots and pseudoknots following below, using the contribution rule:
\[\includegraphics{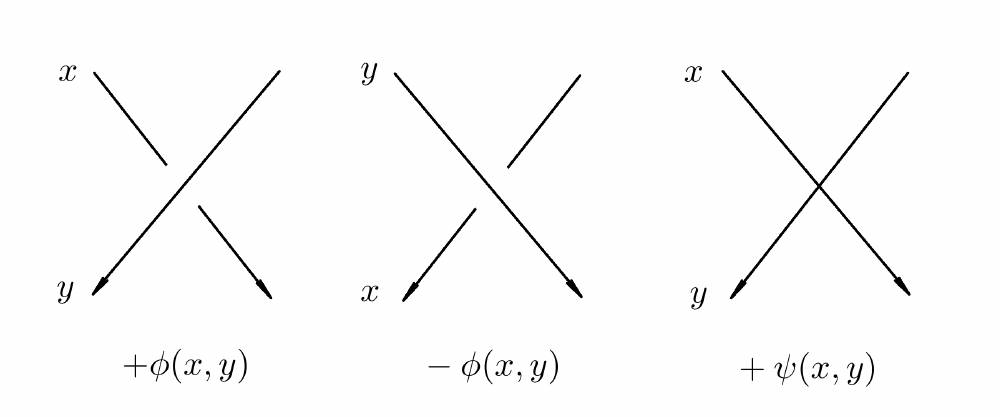}.\]
\[ 
\begin{tikzpicture}[use Hobby shortcut]
	\draw[thick,-stealth] (-6,2) .. (-6,1.8).. (-7,.2) .. (-7,0)..(-7,-.2) .. (-6,-1.8).. (-6,-2);
	\draw [line width=2mm,white,-stealth](-6,0)..(-6,-.2).. (-7,-1.8)..(-7,-2);	
 	\draw[thick,-stealth] (-7,2).. (-7,1.8).. (-6,.2).. (-6,0)..(-6,-.2).. (-7,-1.8)..(-7,-2) ;
	\node[left,lightgray] at (-7,2) {$x$}; 
	\node[left] at (-7,1) {$\psi(x,(y \otr x) \od^{-1}x) $};
	\node[left,lightgray] at (-7,0) {$(y \otr x) \od^{-1}x$};
	\node[left] at (-7,-1) {$\phi((y \otr x) \od^{-1}x,(x \utr y) \od^{-1}y) $};
	\node[left,lightgray] at (-7,-2) {$(x \utr y) \od^{-1}y$};
	\node[right,lightgray] at (-6,2) {$y \otr x$};
	\node[right,lightgray] at (-6,0) {$x \ud ((y \otr x) \otr^{-1}x)$};
	\node[right,lightgray] at (-6,-.5) {$=[(x\utr y)\od^{-1}y]\otr[(y\otr x)\od^{-1}x]$};
	\node[right,lightgray] at (-6,-2) {$y \ud ( x \utr y)$};
	\node[right,lightgray] at (-6.5,-2.5) {$= [(y \otr x) \od^{-1} x] \utr [(x \utr y) \od^{-1} y]$};
	
    \draw [very thick, <->] (-3,0) -- (-1,0);

	\draw[thick,-stealth] (2,2) .. (2,1.8).. (3,.2) .. (3,0)..(3,-.2) .. (2,-1.8).. (2,-2);
	\draw [line width=2mm,white,-stealth](3,2).. (3,1.8).. (2,.2).. (2,0) ;
	\draw[thick,-stealth] (3,2).. (3,1.8).. (2,.2).. (2,0)..(2,-.2).. (3,-1.8)..(3,-2);
	\node[left,lightgray] at (2,2) {$x$}; 
	\node[left] at (2,1) {$\phi(x,y)$};
	\node[left,lightgray] at (2,0) {$y$};
	\node[left] at (2,-1) {$\psi(y,(x \utr y) \od^{-1} y)$};
	\node[left,lightgray] at (2,-2) {$(x \utr y) \od^{-1} y$};
	\node[right,lightgray] at (3,2) {$y \otr x$};
	\node[right,lightgray] at (3,0) {$x \utr y$};
	\node[right,lightgray] at (3,-2) {$y \ud (x \utr y)$};
\end{tikzpicture}\]

\[ 
\begin{tikzpicture}[use Hobby shortcut]
 	\draw[thick,-stealth] (-.5,1).. (-.5,.8) .. (1,0) .. (2,-1) .. (2,-2);
	\draw[thick,-stealth] (1,1) .. (.5,.2).. (-.3,-.5) .. (.6,-1) .. (1,-2);
	\draw [line width=2mm,white,-stealth](2,1) .. (2,0) .. (0,-1) .. (-.5,-2);
	\draw [thick,-stealth](2,1) .. (2,0) .. (0,-1) .. (-.5,-2);

	\node[left,lightgray] at (-.5,1) {$x$};
    \node[left,lightgray] at (-.5,-.5) {$y$};
    \node[left,lightgray] at (-.5,-2) {$z$};
    
    \node[right,lightgray] at (1,1) {$y \od x$};
    \node[right,lightgray] at (1,.1) {$x \ud y$};
    \node[right,lightgray] at (.7,-.8) {$z \utr y$};
    
    \node[left] at (-.3, .2) {$\psi(x,y)$};
    \node[left] at (-.3, -1.1) {$\phi(y,z)$};
    \node[right] at (2,0) {$\phi(x \ud y, z \otr y)$};
    \draw [very thick, <->] (3,-.5) -- (4,-.5);

	\draw[thick,-stealth] (6,1) .. (6,.5).. (7,-.5) .. (6.6,-1) .. (6,-2);   
    \draw[thick,-stealth] (4.5,1).. (4.5,.8) .. (5,0) .. (7,-1.5) .. (7,-2);
    \draw [line width=2mm,white,-stealth] (7,1) .. (7,.5) .. (5,-1) .. (4.5,-2);
    \draw[thick,-stealth] (7,1) .. (7,.5) .. (5,-1) .. (4.5,-2);
    
    \node[left,lightgray] at (4.5,1) {$x$};
    \node[left,lightgray] at (4.5,-2) {$z$};
    
    \node[left,lightgray] at (6,1) {$y \od x$};
    \node[left,lightgray] at (6.4,0) {$z \otr x$};
    \node[left,lightgray] at (6.4,-1) {$x \utr z$};
    \node[left,lightgray] at (6,-2) {$y \utr z$};
    
    \node[left] at (5.5,-.5) {$\phi(x,z)$};
    \node[right] at (7, 0) {$\phi(y \od x, z \otr x)$};
    \node[right] at (7,-1) {$\psi(x \utr z, y \utr z)$};
\end{tikzpicture}\]
\[\begin{tikzpicture}[use Hobby shortcut]
	\draw [thick,-stealth](2,1) .. (2,0) .. (0,-1) .. (-.5,-2);
	\draw [line width=2mm,white,-stealth](-.5,1).. (-.5,.8) .. (1,0) .. (2,-1) .. (2,-2);
	\draw [line width=2mm,white,-stealth](1,1) .. (.5,.2).. (-.3,-.5) .. (.6,-1) .. (1,-2);
 	\draw[thick,-stealth] (-.5,1).. (-.5,.8) .. (1,0) .. (2,-1) .. (2,-2);
	\draw[thick,-stealth] (1,1) .. (.5,.2).. (-.3,-.5) .. (.6,-1) .. (1,-2);
	
	\node[left,lightgray] at (-.5,1) {$z$};
    \node[left,lightgray] at (-.5,-.5) {$y$};
    \node[left,lightgray] at (-.5,-2) {$x$};
    
    \node[right,lightgray] at (1,1) {$y \od z$};
    \node[right,lightgray] at (1,.1) {$z \ud y$};
    \node[right,lightgray] at (.7,-.8) {$x \utr y$};
    
    \node[left] at (-.3, .2) {$\psi(z,y)$};
    \node[left] at (-.3, -1.1) {$-\phi(x,y)$};
    \node[right] at (2,0) {$-\phi(x \utr y, z \ud y)$};
	
    \draw [very thick, <->] (3,-.5) -- (4,-.5);

	\draw[thick,-stealth] (7,1) .. (7,.5) .. (5,-1) .. (4.5,-2);
	\draw [line width=2mm,white,-stealth] (6,1) .. (6,.5).. (7,-.5) .. (6.6,-1) .. (6,-2);
    \draw [line width=2mm,white,-stealth] (4.5,1).. (4.5,.8) .. (5,0) .. (7,-1.5) .. (7,-2);
	\draw[thick,-stealth] (6,1) .. (6,.5).. (7,-.5) .. (6.6,-1) .. (6,-2);   
    \draw[thick,-stealth] (4.5,1).. (4.5,.8) .. (5,0) .. (7,-1.5) .. (7,-2);
    
    \node[left,lightgray] at (4.5,1) {$z$};
    \node[left,lightgray] at (4.5,-2) {$x$};
    
    \node[left,lightgray] at (6,1) {$y \od z$};
    \node[left,lightgray] at (6.5,-.1) {$x \utr z$};
    \node[left,lightgray] at (6.5,-.85) {$z \otr x$};
    \node[left,lightgray] at (6,-2) {$y \otr x$};
    
       \node[left] at (5.5,-.5) {$-\phi(x,z)$};
    \node[right] at (7, 0) {$-\phi(x \utr z, y \od z)$};
    \node[right] at (7,-1) {$\psi(z \otr x, y \otr x)$};
\end{tikzpicture}\]

\[ 
\begin{tikzpicture}[use Hobby shortcut]
\begin{knot}[
consider self intersections=false,
  ignore endpoint intersections=false,
  only when rendering/.style={
  }]
	\strand[thick,-stealth] (-2,2).. (-2,1.8) .. (-2,1) .. (-2,0).. (-1,-1) .. (-.5,0) .. (-1,.5)..(-2,-.5) .. (-2,-2); 
    \end{knot}
    \draw [very thick, <->] (0,0) -- (1,0);
    \begin{knot}[consider self intersections=false,
  ignore endpoint intersections=false,
  only when rendering/.style={
  }]
	\draw[thick,-stealth] (1.5,2) .. (1.5,-2);	
    \end{knot}
    
    \node[left,lightgray] at (-2,2) {$x$};
    \node[right,lightgray] at (-.5,-.5) {$x \ud x=x \od x$};
    \node[left,lightgray] at (-2,-2) {$x$};
    \node[left] at (-2,0) {$\psi(x,x)$};
    \node[right,lightgray] at (1.5,2) {$x$}; 
    
    \begin{knot}[
consider self intersections=false,
  ignore endpoint intersections=false,
  only when rendering/.style={
  }]
  \strand[thick,-stealth] (5,2).. (5,1.8) .. (5,1) .. (5,0).. (4,-1) .. (3.5,0) .. (4,.5)..(5,-.5) .. (5,-2); 
    \end{knot}
    \draw [very thick, <->] (6,0) -- (7,0);
    \begin{knot}[consider self intersections=false,
  ignore endpoint intersections=false,
  only when rendering/.style={
  }]
	\draw[thick,-stealth] (7.5,2) .. (7.5,-2);	
    \end{knot}
    
    \node[left,lightgray] at (5,2) {$x \otr x$};
    \node[left,lightgray] at (3.5,0) {$x$};
    \node[left,lightgray] at (5,-2) {$x \otr x$};
    \node[right] at (5,-.5) {$\psi(x,x)$};
    \node[right,lightgray] at (7.5,2) {$x \otr x$}; 
    
\end{tikzpicture}\]

\begin{definition}
Let $X$ be a psyquandle (respectively, a $pI$-adequate psyquandle), $R$
a commutative ring with identity and $(\phi,\psi)$ a Boltzmann weight 
(respectively, a $pI$-adequate Boltzmann weight). Let $L$ be an oriented 
singular knot or link (respectively, any oriented pseudoknot or pseudolink).
\begin{itemize}
\item[(1)] For each $X$-coloring $L_c$ of $L$ in the set $\mathcal{C}(L,X)$
of $X$-colorings of $L$, we define the \textit{Boltzmann weight} of $L_c$,
denoted $BW(L_c)$, to be the sum over all crossings in $L_c$ 
of crossing contributions as shown:
\[\includegraphics{jc-sn2020-1.pdf}\]
\item[(2)]
We define the \textit{single-variable Boltzmann-enhanced psyquandle polynomial}
to be 
\[\Phi_X^{\phi,\psi}(L)\sum_{L_c\in\mathcal{C}(L,X)} w^{BW(L_c)}\] 
\item[(3)] If $\phi$ and $\psi$ are strongly compatible, we define the 
\textit{partial Boltzmann weights} $BW_\phi(L_c)$ and $BW_\psi(L_c)$ to be 
the sums of $\phi$ contributions and $\psi$ contributions respectively; then
we define the \textit{two-variable Boltzmann-enhanced psyquandle polynomial}
to be
\[\Phi_X^{\phi,\psi}(L)\sum_{L_c\in\mathcal{C}(L,X)} u^{BW_{\phi}(L_c)}v^{BW_{\psi}(L_c)}.\] 
\end{itemize}
\end{definition}

By construction, we have
\begin{proposition}
Let $X$ be a psyquandle, $R$ a commutative ring with identity and
$(\phi,\psi)$ a Boltzmann weight with coefficients in $R$. Then:
\begin{itemize}
\item[(1)] $\Phi_{X}^{\phi,\psi}$ is an invariant of singular knots and links,
\item[(2)] If $X$ and $\psi$ are $pI$-adequate, then 
$\Phi_{X}^{\phi,\psi}$ is an invariant of pseudoknots and pseudolinks.
\end{itemize}
\end{proposition}

\section{Examples}\label{E}

In this section we collect some computations and examples.

\begin{example}\label{ex1}
Let $X=\mathbb{Z}_5$ with the following operations,
\begin{eqnarray*}
x \utr y &=& 3x+4y\\
x \otr y &=& 2x\\
x \ud y &=& 4x+3y\\
x \od y &=& x+y
\end{eqnarray*} 
is an Alexander psyquandle which is \textit{pI-adequate}, for detail see \cite{NOS}. Consider the Boltzmann weight on $X$, defined by $\phi, \psi : X \times X \rightarrow \mathbb{Z}_4$ with $\phi(x,y)=0$ and $\psi(x,y) = 2$ which is \textit{pI-adequate}. The singular knot $K_1$
\[
\begin{tikzpicture}[use Hobby shortcut]
\begin{knot}[
]

\strand[decoration={markings,mark=at position .25 with
    {\arrow[scale=3,>=stealth]{<}}},postaction={decorate}] (-1,0) circle[radius=2cm];
\strand[decoration={markings,mark=at position .25 with
    {\arrow[scale=3,>=stealth]{>}}},postaction={decorate}] (1,0) circle[radius=2cm];
\end{knot}
\node[circle,draw=black, fill=black, inner sep=0pt,minimum size=6pt] (a) at (0,1.7) {};
\node[circle,draw=black, fill=black, inner sep=0pt,minimum size=6pt] (a) at (0,-1.7) {};
\node[left] at (-3,0) {$x_1$};
\node[left] at (-1,0) {$x_2$};
\node[right] at (1,0) {$x_3$};
\node[right] at (3,0) {$x_4$};
\end{tikzpicture}
\]
has the system of coloring equations given by 
\begin{eqnarray*}
4x_1 -2x_2=x_1 \ud x_2 &=& x_3\\
x_1+x_2 = x_2 \od x_1 &=& x_4\\
x_1+x_2 = x_1 \od x_2 &=& x_3\\
-2x_1+4x_2 = x_2 \ud x_1 &=& x_4
\end{eqnarray*}
which we can solve by row-reduction over $\mathbb{Z}_5$:
\[ 
\left[
\begin{array}{cccc}
 4 & -2 & -1 & 0 \\
 1 & 1 & 0 & -1 \\
 1 & 1 & -1 & 0 \\
 -2 & 4 & 0 & -1 \\
\end{array}
\right]
\rightarrow
\left[
\begin{array}{cccc}
 1 & 0 & 0 & 2 \\
 0 & 1 & 0 & 2 \\
 0 & 0 & 1 & 4 \\
 0 & 0 & 0 & 0 \\
\end{array}
\right].
\]

Now, consider the following singular knot $K_2$ 
\[
\begin{tikzpicture}[use Hobby shortcut]
\begin{knot}[
]

\strand[decoration={markings,mark=at position .25 with
    {\arrow[scale=3,>=stealth]{<}}},postaction={decorate}] (-1,0) circle[radius=2cm];
\strand[decoration={markings,mark=at position .25 with
    {\arrow[scale=3,>=stealth]{>}}},postaction={decorate}] (1,0) circle[radius=2cm];
\end{knot}
\node[circle,draw=black, fill=black, inner sep=0pt,minimum size=6pt] (a) at (0,1.7) {};

\node[left] at (-3,0) {$x_1$};
\node[left] at (-1,0) {$x_2$};
\node[right] at (1,0) {$x_3$};
\node[right] at (3,0) {$x_4$};
\end{tikzpicture}
\]
has the system of coloring equations given by 
\begin{eqnarray*}
4x_1 -2x_2=x_1 \ud x_2 &=& x_3\\
x_1+x_2 = x_2 \od x_1 &=& x_4\\
2x_1 = x_1 \otr x_2 &=& x_3\\
-x_1+3x_2 = x_2 \utr x_1 &=& x_4
\end{eqnarray*}
which we can solve by row-reduction over $\mathbb{Z}_5$:
\[
\left[
\begin{array}{cccc}
 4 & -2 & -1 & 0 \\
 1 & 1 & 0 & -1 \\
 2 & 0 & -1 & 0 \\
 -1 & 3 & 0 & -1 \\
\end{array}
\right]
\rightarrow
\left[
\begin{array}{cccc}
 1 & 0 & 0 & 2 \\
 0 & 1 & 0 & 2 \\
 0 & 0 & 1 & 4 \\
 0 & 0 & 0 & 0 \\
\end{array}
\right].
\]
These two systems have 5 solutions, therefore, both diagrams have 5 coloring 
by this psyquandle. We will now consider the Boltzmann weight enhanced 
invariant in order to distinguish these two singular knots. Using the
Boltzmann weight above, we obtain enhanced invariant values
$\Phi_{X}^{\phi,\psi}(K_1)  = 5$ and $\Phi_{X}^{\phi,\psi}(K_2)= 5v^2$, demonstrating 
that the enhanced invariant is not determined by the number of colorings
and hence is a proper enhancement. 
\end{example}

\begin{example}\label{ex2}
Using our custom \texttt{Python} code, we computed the Boltzmann weight for certain singular knots known as \emph{two-bouquet graphs of type $K$} (with choice of orientation) listed in \cite{O} using the psyquandle $X$ given by the operation matrix 
\[\left[\begin{array}{cc|cc|cc|cc}
1 & 1 & 1 & 1 &  2 & 2 & 2 & 2 \\ 
2 & 2 & 2 & 2 &  1 & 1 & 1 & 1

\end{array}\right].\]
The weight function $\phi: X \times X \rightarrow \mathbb{Z}_{14}$ is given by the following matrix
\[\left[\begin{array}{cccc}
0 & 0  \\ 
7 & 0  
\end{array}\right]\]
and the weight function $\psi: X \times X \rightarrow \mathbb{Z}_{14}$ is given by the following matrix
\[\left[\begin{array}{cccc}
0 & 0 \\ 
3 & 0 
\end{array}\right].\]

We compute the Boltzmann weight enhancement for $5_2^k$ and $5_3^k$ from \cite{O}. We obtain the following Boltzmann weight enhanced invariant values: $\Phi_X^{\phi,\psi}(5_2^k)  = 2u^7$ and $\Phi_X^{\phi,\psi}(5_3^k)  = 2$. This example shows that $\Phi_X^{\phi,\psi}$ detects additional information beyond counting invariant $\Phi_X^\mathbb{Z}(5_2^k) = 2 = \Phi_X^\mathbb{Z}(5_3^k)$. 
\end{example}

\begin{example}\label{exlist}
Using our custom \texttt{Python} code, we compute the Boltzmann weights invariant for the \emph{2-bouquet graphs of type $L$} (with choice of orientation) in \cite{O} using the psyquandle given by the operation matrix
\[\left[\begin{array}{cccccc|cccccc|cccccc|cccccc}
2& 4 & 4 & 6 & 6 & 2 &  2 & 6 & 2 & 6 & 2 & 6 &  2& 4& 2& 6& 2& 2 &  2& 6& 4& 6& 6& 6\\ 
3& 5 & 5 & 1 & 1 & 3 &  1 & 5 & 1 & 5 & 1 & 5 &  3& 5& 5& 5& 1& 5 &  1& 5& 1& 1& 1& 3\\
4& 6 & 6 & 2 & 2 & 4 &  6 & 4 & 6 & 4 & 6 & 4 &  6& 6& 6& 2& 6& 4 &  4& 4& 6& 4& 2& 4\\
5& 1 & 1 & 3 & 3 & 5 &  5 & 3 & 5 & 3 & 5 & 3 &  5& 3& 1& 3& 3& 3 &  5& 1& 5& 3& 5& 5\\
6& 2 & 2 & 4 & 4 & 6 &  4 & 2 & 4 & 2 & 4 & 2 &  4& 2& 4& 4& 4& 6 &  6& 2& 2& 2& 4& 2\\
1& 3 & 3 & 5 & 5 & 1 &  3 & 1 & 3 & 1 & 3 & 1 &  1& 1& 3& 1& 5& 1 &  3& 3& 3& 5& 3& 1
\end{array}\right]\]
with weight function $\phi: X \times X \rightarrow \mathbb{Z}_2$ given by the matrix
\[\left[\begin{array}{cccccc}
0& 1& 0& 1& 0& 1\\
0& 0& 0& 0& 0& 0\\
0& 1& 0& 1& 0& 1\\
0& 0& 0& 0& 0& 0\\
0& 1& 0& 1& 0& 1\\
0& 0& 0& 0& 0& 0
\end{array}\right]\]
and weight function $\psi: X \times X \rightarrow \mathbb{Z}_2$ given by the matrix
\[\left[\begin{array}{cccccc}
1& 0& 1& 0& 1& 0\\
1& 1& 1& 1& 1& 1\\
1& 0& 1& 0& 1& 0\\
1& 1& 1& 1& 1& 1\\
1& 0& 1& 0& 1& 0\\
1& 1& 1& 1& 1& 1
\end{array}\right].\]
The results are collected in the table

\begin{center}
\begin{tabular}{ r| r | l }
$\Phi_X^\mathbb{Z}$ & $\Phi_X^{\phi,\psi}(L)$ & $L$ \\
\hline
     $12$	& $12w$ 	& 	$5_2^l,6_1^l$ \\
  		& $6w+6$	& 	$3_1^l,4_1^l,5_3^l,6_2^l,6_6^l$\\
\hline
     $24$ 	& $24w$ 	&    	$6_3^l,6_8^l,6_9^l,6_{10}^l,6_{11}^l$\\	
 	    	& $6w+18$ 	&   	$5_1^l,6_5^l,6_7^l$ \\
  		& $18w+6$ 	&  	$1_1^l$\\
\hline
	$36$ & $18w+18$ & 	$6_4^l, 6_{12}^l$.
\end{tabular}
\end{center}
\end{example}
\begin{example}\label{plistex}
Using our custom \texttt{Python} code, we computed the counting invariant 
$\Phi_X^\mathbb{Z}$ and Boltzmann weight enhanced invariant $\Phi_X^{\phi,\psi}$ 
for a choice of orientation for the pseudoknots in \cite{HHJJMR} using the 
psyquandle $X$ given by the operation matrix 
\[\left[\begin{array}{ccc|ccc|ccc|ccc}
1& 3& 1  &  1& 1& 1  &  3& 1& 3  &  3& 3& 3\\ 
2& 2& 2  &  2& 2& 2  &  2& 2& 2  &  2& 2& 2\\
3& 1& 3  &  3& 3& 3  &  1& 3& 1  &  1& 1& 1
\end{array}\right]\]
and the weight function $\phi: X \times X \rightarrow \mathbb{Z}_6$ given by the matrix
\[\left[\begin{array}{ccc}
0& 3& 0\\
0& 0& 0\\
0& 0& 0
\end{array}\right]\]
and weight function $\psi: X \times X \rightarrow \mathbb{Z}_6$ given by the matrix
\[\left[\begin{array}{ccc}
0& 5& 4\\
2& 0& 5\\
4& 5& 0
\end{array}\right].\]
We have a pI-adaquate psyquandle since the two right blocks have the same diagonal and the pair $(\phi,\psi)$ is pI-adequate since the entries along the diagonal of $\psi$ are all zero.
The results are collected in the table

\begin{center}
\begin{tabular}{ r| r | l }
$\Phi_X^\mathbb{Z}$ & $\Phi_X^{\phi,\psi}(L)$ & $L$ \\
 	\hline
  $3$	& $3$ 		& $3_1.1,3_1.3,4_1.1,4_1.3,4_1.5,5_1.1,5_1.3,5_1.5,5_2.1,5_2.3,5_2.5,5_2.6,5_2.8,5_2.10$ \\
  		& $2w^2+1$	& $3_1.2,4_1.2,4_1.4,5_1.4,5_2.4,5_2.7,5_2.9 $\\
  		& $2w^4+1$ 	& $5_1.2,5_2.2$.
\end{tabular}
\end{center}

\begin{example}\label{exjablan}
Using the psyquandle and weight functions from Example~\ref{plistex} the Boltzmann weight of a pseudoknot can detect additional information for pseudoknots with the same Jablan polynomial. We collect the results in the following two tables
\begin{center}
\begin{tabular}{ r| r | l }
$\Delta_J(L)$ & $\Phi_X^{\phi,\psi}(L)$ & $L$ \\
 	\hline
  $1$	& $3$ 		& $3_1.1, 4_1.1, 4_1.3,5_1.1,5_2.1,5_2.6$ \\
  		& $2w^2+1$	& $3_1.2, 4_1.2,5_2.9$ \\
  		& $2w^4+1$ 	& $5_1.2,5_2.2$  		
\end{tabular}
\end{center}
and
\begin{center}
\begin{tabular}{ r| r | l }
$\Delta_J(L)$ & $\Phi_X^{\phi,\psi}(L)$ & $L$ \\
 	\hline
  $s^2+t^2$	& $3$ 		& $3_1.3, 5_2.3$ \\
  		& $2w^2+1$	& $5_2.4.$  		
\end{tabular}
\end{center}
\end{example}
\end{example}

\section{Questions}\label{Q}

We end with a few questions for future work.

In \cite{CN1}, biquandle homology was enhanced for virtual knots and links
with a weight function at virtual crossings, and these were found to 
be cocycles themselves in a cohomology theory we called \textit{$S$-cohomology}.
Is something similar true for psyquandles? How should the singular crossing 
weight $\psi$ be interpreted in terms of cohomology?

What other enhancements of the psyquandle counting invariant and 
Boltzmann weight-enhanced invariants cane be defined?

\bibliography{jc-sn2020}{}
\bibliographystyle{abbrv}

\bigskip

\noindent
\textsc{
}

\bigskip

\noindent
\textsc{Department of Mathematical Sciences \\
Claremont McKenna College \\
850 Columbia Ave. \\
Claremont, CA 91711}\\\\
\textsc{Mathematics and Statistics Department \\
Hamilton  College \\
198 College Hill Rd. \\
Clinton, NY 13323}
\end{document}